\documentclass[12pt,a4paper]{article}
\usepackage{eurosym}
\usepackage{amssymb}
\usepackage{amsmath}

\newtheorem{theorem}{Theorem}

\newtheorem{condition}[theorem]{Condition}

\newtheorem{corollary}[theorem]{Corollary}

\newtheorem{definition}[theorem]{Definition}

\newtheorem{lemma}[theorem]{Lemma}

\newtheorem{proposition}[theorem]{Proposition}
\newtheorem{remark}[theorem]{Remark}

\numberwithin{equation}{section}

\begin{document}

\title{Stochastic Camassa-Holm equation with convection type noise \footnote{{\bf Mathematics Subject Classification (2010)}. Primary 60H15; Secondary 60H25, 35R60, 76B15, 35Q86.
{\bf Keywords.} Stochastic  Camassa-Holm equation, local strong solutions, random quasi-linear partial differential equations, Kato's operator theory method, Doss-Sussman correspondence.}
}
\author{Sergio Albeverio \\
{\small Institute for Applied Mathematics, }\\
{\small Rheinische Friedrich-Wilhelms Universit\"{a}t Bonn, }\\
{\small and Hausdorff Center for Mathematics, }\\
{\small Endenicher Allee 60, D-53115 Bonn, Germany} \and Zdzis\l aw Brze\'{z}%
niak, Alexei Daletskii \\
{\small Department of Mathematics, University of York, }\\
{\small Heslington, York, YO10 5DD, UK }
}


\maketitle

\begin{abstract}
We consider a stochastic Camassa-Holm equation driven by a one-dimensional
Wiener process with a first order differential operator as diffusion
coefficient. We prove the existence and uniqueness of local strong solutions
of this equation. In order to do so, we transform it into a random
quasi-linear partial differential equation and apply Kato's operator theory
methods. Some of the results have potential to find applications to other
nonlinear stochastic partial differential equations.
\end{abstract}

\section{Introduction}

The (deterministic) Camassa-Holm (CH) equation is a non-local partial
differential equation describing propagation of waves in shallow water.
Although fist introduced by B. Fuchssteiner and A. Fokas in \cite{FF} as
part of a family of integrable Hamiltonian equations, it was rediscovered by
R. Camassa and D. Holm \cite{CaHo}, who gave its physical derivation and
interpretation. In contrast to the Korteveg-de-Vries equation, the CH
equation admits so-called peaked solutions describing wave breaking
phenomena. Various aspects of the CH equation have been extensively studied,
see literature reviews in, e.g., \cite{CE} and \cite{BHP}. In particular, it
is known that the CH equation is locally well-posed in Sobolev spaces $H^{s}$%
, $s>3/2$ (here and in what follows, $H^{n}:=H^{n,2}(\mathbb{R})$, $n\in
\mathbb{N},$ is the real Sobolev space). Depending on the shape of the
initial data, the solution can either exist globally or blow up in any
Sobolev space, with its slope becoming vertical in finite time \cite{CE}.

The CH equation has the form%
\begin{multline}
u_{t}-u_{xxt}+3uu_{x}-2u_{x}u_{xx}-uu_{xxx}  \label{CHeq} \\
\equiv (1-\partial _{x}^{2})u_{t}+(1-\partial _{x}^{2})\left[ uu_{x}\right]
+\partial _{x}\left[ u^{2}+\frac{1}{2}\left( u_{x}\right) ^{2}\right] =0,\
t>0,x\in \mathbb{R},
\end{multline}%
where $u(t,x)$ denotes the fluid velocity at time $t$. Here $\partial _{x}:=%
\frac{\partial }{\partial x}$. Introducing a momentum density
\begin{equation*}
y:=u-u_{xx}\equiv (1-\partial _{x}^{2})u=:Q^{2}u,
\end{equation*}%
one can rewrite equation (\ref{CHeq}) in a quasi-linear form%
\begin{equation}
y_{t}(t)+A(y(t))y(t)=0  \label{QLE}
\end{equation}%
in $L^{2}:=L^{2}(\mathbb{R})$ (or any other suitable functional space). Here
$A(v):=a(v)\partial _{x}+b(v)$, $v\in H^{1}$, is the first-order
differential operator with coefficients $a(v)=Q^{-2}v\in H^{3},\
b(v)=2\left( \partial _{x}Q^{-2}v\right) \in H^{2}$, that is,%
\begin{equation}
\left[ A(v)f\right] (x)=a(v)\partial _{x}f(x)+b(v)f(x),\ x\in \mathbb{R},\
f\in H^{1}.  \label{A-CH}
\end{equation}

Recently, D. Holm \cite{Holm} proposed an approach for including stochastic
perturbations in hydrodynamics equations. This approach is based on a
stochastic extension of the variational principle in fluid dynamics. The
corresponding stochastic version of the CH equation (\ref{CHeq}) was
introduced D. Crisan and D. Holm in \cite{CrHo}. It has the following form:
\begin{multline}
dy(t)+A(y(t))y(t)dt+\sum_{k=1}^{n}\left( \partial _{x}y(t)+y(t)\partial
_{x}\right) \xi _{k}\circ dw_{k}(t)  \label{SCH} \\
\equiv dy(t)+A(y(t))y(t)dt+\sum_{k=1}^{n}D_{k}y(t)\circ dw_{k}(t)=0,\
t>0,y\in \mathbb{R}.
\end{multline}%
Here $D_{k}:=\xi _{k}\partial _{x}+\partial _{x}\xi _{k}\equiv \xi
_{k}\partial _{x}+2\left( \partial _{x}\xi _{k}\right) \mathrm{,}$ $k=1,..n,$
are first-order differential operators associated with suitable functions
(vector fields) $\xi _{k}:\mathbb{R}\rightarrow \mathbb{R}$, $w_{k}$, $%
k=1,..n,$ are independent Wiener processes and $\circ dw_{k}(t)$ stands for
the Stratonovich stochastic differential (see Def.~\ref{strong-sol} below).
For further developments from \cite{Holm}, \cite{CrHo} see, e.g., \cite%
{Holm1}, \cite{Holm-Ivanov}; these also relate to stochastic thermodynamics
and turbulence, for which we refer to e.g. \cite{BCF1}, \cite{CFHOTW}, \cite%
{E}, \cite{Fland}.

In this paper, we study the case of a single vector field $\xi $ (with $n=1$%
). In order to deal with the diffusion term of equation (\ref{SCH}), we
transform it into a partial differential equation with random coefficients.
This approach goes back to the paper \cite{Sus} by Sussman , see also Doss
\cite{Doss}. These works were concerned with stochastic ordinary
differential equations and motivated by the control theory. In stochastic
partial differential equations (SPDEs) theory, the Doss-Sussman method was
first used in \cite{101} and \cite{BCF}. Both papers studied the Wong-Zakai
approximations (or robustness) of linear SPDEs with drift being the
generator of an analytic semigroup. The corresponding Banach space setting
generalizations can be found in \cite{103}.

Recently, the Doss-Sussman method was used in \cite{GLT} to study the
convergence of a finite element method for stochastic
Landau-Lifshitz-Gilbert equations. The Wong-Zakai approximations to such
equations were studied in \cite{BMM}. Other related papers are \cite{104}
and \cite{102}, where it was noted that the Doss-Sussman method could lead
to an alternative proof of the main result therein, and then applied to
(nonlinear) stochastic compressible Euler equations, respectively. Another
example of the use of the Doss-Sussman is its application to the stochastic
nonlinear Schr{\"{o}}dinger equation, see \cite{104}.

After applying the Doss-Sussman method to equation (\ref{SCH}) (with $n=1$),
we study the resulting partial differential equation (PDE) using a modified
version of the approach of \cite{CE} based on the Kato operator theory
techniques. Our main result is the proof of the existence and uniqueness of
local strong solutions of equation (\ref{SCH}). We hope that with an
additional work it should also be possible to prove the robustness of this
equation. Also, a modification of our method should allow for the study of
the case of multiple (non-commuting) vector fields $\xi _{k}$ on the
right-hand side of (\ref{SCH}).

Let us mention that a stochastic CH equation with additive noise was
introduced and studied in \cite{CGG}; the case of a multiplicative noise
given by a one-dimensional Wiener process with $H^{s}$-continuous diffusion
coefficient was considered in \cite{CG} and \cite{Tan}. Those studies do not
cover the case of the noise as in (\ref{SCH}), where the diffusion
coefficient is generated by an unbounded linear operator. The importance of
studying equation (\ref{SCH}) has been stressed by D. Crisan and D. Holm in
\cite{CrHo} because of its geometric and physical motivations, and its
relevance in geophysical applications.

The structure of the paper is as follows. In Section \ref{sec1}, we
formulate the main result and derive the explicit form of the PDE obtained
by the Doss-Sussman method. Section \ref{sec2} is devoted to the general
Kato method and its application to the latter PDE, which leads to the proof
of our main result in Section \ref{main-proof}. In Section \ref{sog} we
provide the proofs of (auxiliary) technical results on the regularity of
one-parameter groups generated by first order differential operators.
Finally, in Section \ref{S-O} we prove the Doss-Sussman correspondence
between SDEs and (random) ordinary differential equations in Hilbert spaces,
adapted to our setting.

\textbf{Acknowledgement. }We are very gratefull to Darryl Holm for his
interest to this work and stimulating discussions.  Part of this research was
carried out during AD's stay at Mathematical Institute of the University of Bonn. Financial support of this stay by Alexander von
Humboldt Foundation is gratefully acknowledged 

\section{Stochastic Camassa-Holm equation\label{sec1}}

\subsection{Formulation of the main result\label{main-result}}

We will consider a stochastic Camassa-Holm equation (SCH) of the form%
\begin{equation}
dy(t)+F(y(t))dt+Dy(t)\circ dw(t)=0,\ t\geq 0,\ F(y):=A(y)v,v\in H^{1},
\label{SPDE}
\end{equation}%
on a suitable filtered probability space $\left( \Omega ,\mathcal{F},\left(
\mathcal{F}_{t}\right) _{t\geq 0},\mathbb{P}\right) $, where $A(v)$ is given
by formula (\ref{A-CH}), $D=\xi \partial _{x}+\eta $, $\xi \in C_{b}^{4}$, $%
\eta \in C_{b}^{3}$ and $w$ is a one-dimensional Wiener process. We will be
looking for a solution of this equation in $H^{2}$.

\begin{definition}
\label{strong-sol}A strong solution of equation (\ref{SPDE}) is an $H^{2}$%
-valued continuous process $y(t)$, $t\in \left[ 0,\theta \right] $, where $%
\theta $ is a finite stopping time, such that the equality%
\begin{equation*}
y(t\wedge \theta )=y_{0}+\int_{0}^{t\wedge \theta
}F(y(s))ds+\int_{0}^{t\wedge \theta }Dy(s)\circ dw(s),\ t\geq 0,
\end{equation*}%
is satisfied in $L^{2}$, $\mathbb{P}$-a.s., for every $t\geq 0$. Here $\circ
dw(s)$ stands for the Stratonovich stochastic differential, that is,
\begin{equation*}
\int_{0}^{t\wedge \theta }Dy(s)\circ dw(s)=\frac{1}{2}\int_{0}^{t\wedge
\theta }D^{2}y(s)ds+\int_{0}^{t\wedge \theta }Dy(s)dw(s).
\end{equation*}
\end{definition}

We can now formulate the main result of this work.

\begin{theorem}
\label{theor-main}For any $y_{0}\in H^{2}$ there exists a stopping time $%
\theta >0$ and a strong solution $y(t)\in H^{2}$, $t\in \left[ 0,\theta %
\right] $, of equation (\ref{SPDE}) with initial condition $y(0)=y_{0}$. If $%
y^{1}(t)$ and $y^{2}(t)$, $t\in \left[ 0,\theta \right] $, are two such
solutions then $y^{1}=y^{2}$.
\end{theorem}

The proof will go along the following lines: first, we reduce equation (\ref%
{SPDE}) to a PDE of a form similar to (\ref{QLE}) but with time-dependent
coefficients, and then apply the general Kato method, similar to the case of
the deterministic Camassa-Holm equation.

\subsection{Reduction to a random PDE\label{diff-pert}}

Let us fix $\xi \in C_{b}^{4},\ \eta \in C_{b}^{3}$ and consider the
one-parametric group $U=\left( U_{t}^{\xi ,\eta }\right) _{t\in \mathbb{R}}$
of operators in $L^{2}$ defined by the formula
\begin{equation}
\left[ U_{t}f\right] (x)=e^{c(t,x)}f(\varphi _{-t}(x)),\ f\in L^{2},\ x\in
\mathbb{R},\ t\geq 0,  \label{U-explicit}
\end{equation}%
where $\varphi _{t}$ is a diffeomorphism generated by the vector field $\xi
\partial _{x}$ and $c(t,x)=\int_{0}^{t}\eta (\varphi _{s-t}(x))ds$, see
Lemma \ref{U} in Section \ref{sog} below.

According to the results of Section \ref{sog}, $U$ is strongly continuous in
the Hilbert spaces $X:=H^{1}$ and $Y:=H^{2}$. For the corresponding
generators $(D^{X},Dom(D^{X}))$ and $(D^{Y},Dom(D^{Y}))$ we have
\begin{equation*}
H^{2}\subset Dom(D^{X})\text{ and }H^{3}\subset Dom(D^{Y}),
\end{equation*}%
and the restrictions of $D^{X}$ and $D^{Y}$ on $H^{2}$ and $H^{3}$,
respectively, coincide with the first order differential operator $D=\xi
\partial _{x}+\eta $. Note also that $Dom(D^{Y})\subset Dom(D^{X})$. It is
shown in Lemma \ref{semi-lemma} below that $U$ satisfies the estimate%
\begin{equation}
\left\Vert U_{t}^{\xi ,\eta }\right\Vert _{\mathcal{L}(X)},\left\Vert
U_{t}^{\xi ,\eta }\right\Vert _{\mathcal{L}(Y)}\leq C_{1}e^{C_{2}\left\vert
t\right\vert },\ t\in \mathbb{R},  \label{U-est0}
\end{equation}%
for some constants $C_{1},C_{2}<\infty $. In this section, we will write $%
U_{t}$ in place of $U_{t}^{\xi ,\eta }$, whenever possible.

According to the results of Section \ref{S-O} (with $Y=H^{2}$ and $\mathfrak{%
X}=L^{2}$), equation (\ref{SPDE}) is equivalent to the following random
integral equation in $L^{2}$:%
\begin{equation}
z(t)=z(0)-\int_{0}^{t}\widehat{F}(s,z(s))ds,\ t\geq 0,  \label{SCH-int}
\end{equation}%
where
\begin{equation*}
\widehat{F}(t,z):=U_{w(t)}F\left( U_{w(t)}^{-1}z\right) \equiv \widehat{A}%
(w(t),z)z,\ t\geq 0,z\in H^{2},
\end{equation*}%
and%
\begin{equation}
\widehat{A}(t,v):=U_{t}A(U_{t}^{-1}v)U_{t}^{-1},\ t\geq 0,v\in H^{2}.
\label{Atv}
\end{equation}%
Our next goal is to study the structure of operator $\widehat{A}(t,v)$.

Consider a generic first order differential operator $\mathcal{A}%
=a_{0}\partial _{x}+b_{0}$ with the coefficients $a_{0}\in H^{3}$ and $%
b_{0}\in H^{2}$ and define operators
\begin{equation}
C(t):=U_{t}\mathcal{A}U_{t}^{-1},\ t\geq 0,  \label{Ct-op}
\end{equation}%
on the domain $H^{2}$. Note that $C(t)\in \mathcal{L}(H^{2},H^{1})$.

\begin{lemma}
\label{C-lemma}Assume that $a_{0}\in H^{3}$, $b_{0}\in H^{2}$. Then operator
$C(t)$ defined above by formula (\ref{Ct-op}) has the form%
\begin{equation}
C(t)v=a(t,\cdot )\partial _{x}v+b(t,\cdot )v,\ v\in H^{2},\ t\geq 0,
\label{C-OD}
\end{equation}%
where $a(t,x)$ and $b(t,x)$ solve the system of first order partial
differential equations%
\begin{eqnarray*}
a_{t}(t,x) &=&\xi (x)a_{x}(t,x)-\xi _{x}(x)a(t,x),\ a(0,x)=a_{0}(x), \\
b_{t}(t,x) &=&\xi (x)b_{x}(t,x)-\eta _{x}(x)a(t,x),\ b(0,x)=b_{0}(x),
\end{eqnarray*}%
(with subscript $x$ denoting as usual the derivative $\partial _{x}$).
Moreover, $a(t):=a(t,\cdot )\in H^{3}$ and $b(t):=b(t,\cdot )\in H^{2}$ and
\begin{equation}
\left\Vert a(t)\right\Vert _{H^{3}}\leq C_{1}e^{tC_{2}}\left\Vert
a_{0}\right\Vert _{H^{3}},\ \left\Vert b(t)\right\Vert _{H^{2}}\leq
C_{1}e^{tC_{2}}\left( \left\Vert a_{0}\right\Vert _{H^{3}}+\left\Vert
b_{0}\right\Vert _{H^{2}}\right)  \label{ab-bound}
\end{equation}%
for some constants $C_{1},C_{2}>0$ (depending only on $\xi $ and $\eta $).
\end{lemma}

\noindent \textbf{Proof.} Let us fix $f\in H^{3}$ and consider the map $%
\mathbb{R}\ni t\mapsto C(t)f\in H^{1}$. Observe that $H^{3}\mathcal{\subset }%
Dom_{Y}(D)$ and $\mathcal{A}U_{t}^{-1}:H^{3}\mathcal{\rightarrow }Dom_{X}(D)$
for all $t$, which implies that the function $C(t)f$, $t\geq 0$, with $C(t)\
$being given by (\ref{Ct-op}), is differentiable and satisfies equation%
\begin{equation*}
\frac{d}{dt}C(t)f=\left[ D,C(t)\right] f,
\end{equation*}%
where $\left[ \cdot ,\cdot \right] $ stands for the commutator. The
substitution of the explicit expression (\ref{U-explicit}) in the formula $%
C(t)=U_{t}\mathcal{A}U_{t}^{-1}$ shows that $C(t)$ has the form (\ref{C-OD}%
). We can now compute the commutator:%
\begin{equation*}
\left[ D,C(t)\right] =\left[ \xi \partial _{x}+\eta ,a(t)\partial _{x}+b(t)%
\right] =\alpha (t)\partial _{x}+\beta (t),
\end{equation*}%
where%
\begin{equation*}
\alpha (t)=\xi a_{x}(t)-\xi _{x}a(t),\ \beta (t)=\xi b_{x}(t)-\eta _{x}a(t).
\end{equation*}%
Observe that $f\in H^{3}$ belongs to the domain of the operators $DC(t)$ and
$C(t)D$. Thus we have%
\begin{equation*}
\frac{d}{dt}C(t)f=\left( \xi a_{x}(t)-\xi _{x}a(t)\right) \partial
_{x}f+\left( \xi b_{x}(t)-\eta _{x}a(t)\right) f,\ t\geq 0,
\end{equation*}%
On the other hand, by (\ref{C-OD}),
\begin{equation*}
\frac{d}{dt}C(t)f=a_{t}\partial _{x}f+b_{t}f,\ t\geq 0,
\end{equation*}%
so that%
\begin{eqnarray*}
a_{t} &=&\xi a_{x}-\xi _{x}a,\ a(0)=a_{0},\ t\geq 0, \\
b_{t} &=&\xi b_{x}-\eta _{x}a,\ b(0)=b_{0},\ t\geq 0.
\end{eqnarray*}%
Thus for any $t\geq 0$ we have
\begin{equation}
C(t)f=a(t)\partial _{x}f+b(t)f,\ f\in H^{3}.  \label{C-OP1}
\end{equation}%
Observe, on the other hand, that the operators $C(t)$ and $a(t)\partial
_{x}+b(t)$ belong to $\mathcal{L}(H^{2},H^{1})$. Thus equality (\ref{C-OP1})
can be extended to any $f\in H^{2}$.

Thus, recalling that $U_{t}^{\xi ,-\xi ^{\prime }}$ is a one-parameter group
generated by the operator $\xi \partial _{x}-\xi _{x}$, we have the
representation%
\begin{equation}
a(t)=U_{t}^{\xi ,-\xi ^{\prime }}a_{0},\ t\geq 0,  \label{at}
\end{equation}%
and
\begin{equation}
b(t)=U_{t}^{\xi ,0}b_{0}+\int_{0}^{t}U_{t-\tau }^{\xi ,0}\left( \eta
^{\prime }a(\tau )\right) d\tau ,\ t\geq 0.  \label{bt}
\end{equation}%
Since by Lemma \ref{semi-lemma} below both $U_{t}^{\xi ,-\xi ^{\prime }}$
and $U_{t}^{\xi ,0}$ leave the spaces $H^{1},H^{2}$ and $H^{3}$ invariant,
we infer that $a(t)\in H^{3}$ and $b(t)\in H^{2}$. The bound (\ref{ab-bound}%
) follows now easily from (\ref{U-est0}), (\ref{at}) and (\ref{bt}). The
proof is complete.

\hfill%
$\square $

We can now return to the operator family $\widehat{A}(t,v)$ given by (\ref%
{Atv}).

\begin{proposition}
\label{ab-pert}For any $t\geq 0$ and $v\in H^{1}$, the operator $\widehat{A}%
(t,v)$ has the form
\begin{equation*}
\widehat{A}(t,v)=a(t,v)\partial _{x}+b(t,v)
\end{equation*}%
on the domain $H^{2}$, where $a(t,v)\in H^{3}$ and $b(t,v)\in H^{2}$ are
given by formulae (\ref{at}) and (\ref{bt}) with
\begin{equation}
a_{0}=Q^{-2}U_{t}^{-1}v\ \text{and }b_{0}=2\partial _{x}Q^{-2}U_{t}^{-1}v,
\label{ab-init}
\end{equation}%
respectively, and satisfy the bound%
\begin{equation}
\left\Vert a(t)\right\Vert _{H^{3}}\leq C_{1}e^{tC_{2}}\left\Vert
v\right\Vert _{H^{1}},\ \left\Vert b(t)\right\Vert _{H^{2}}\leq
C_{1}e^{tC_{2}}\left\Vert v\right\Vert _{H^{1}}.  \label{ab-v-bound}
\end{equation}%
for some $C_{1},C_{2}<\infty $.
\end{proposition}

\noindent \textbf{Proof.} We can first fix any $s$ and apply Lemma \ref%
{C-lemma} to operator (\ref{Ct-op}) with $\mathcal{A}:\mathcal{=}%
A(U_{s}^{-1}v)$ and then set $s=t$. The bound (\ref{ab-v-bound}) follows
from (\ref{ab-bound}) and estimate (\ref{U-est}) of the norm of $U_{t}$.
\hfill%
$\square $

\begin{corollary}
\label{a-cont}For any $v\in H^{1}$ we have $\widehat{A}(t,v)\in \mathcal{L}%
(H^{2},H^{1})$ and the map $\mathbb{R}\ni t\mapsto \widehat{A}(t,v)\in
\mathcal{L}(H^{2},H^{1})$ is continuous.
\end{corollary}

\noindent \textbf{Proof. }The result follows from formulae (\ref{at}), (\ref%
{bt}), (\ref{ab-init}) and the strong continuity of the one-parameter groups
$U_{t}^{\xi ,-\xi ^{\prime }}$, $U_{t}^{\xi ,0}$ and $U_{t}$.
\hfill%
$\square $

\section{Quasi-linear equations via Kato's method\label{sec2}}

\subsection{General Kato's method}

Consider a pair of densely embedded Hilbert spaces $Y\subset X$ and a
quasi-linear equation in $X$:%
\begin{equation}
\frac{d}{dt}v+A(t,v)v=0,\ v(0)=v_{0}\in Y,\ t\in \left[ 0,T\right] ,
\label{ODE_semi}
\end{equation}%
for some $T>0$, where $A(t,v)$ is a linear (unbounded) operator in $X$ with
domain $D_{t,v}:=Dom(A(t,v))\supset Y$.

We introduce the following condition, which is a version of the condition
given in \cite[page 34]{Kato} adapted to our setting . Let $I\subset \mathbb{%
R}$\textrm{\ }be an interval.

\begin{condition}
\label{Kato-cond} There exists $R>0$ such that the operator family $A(t,v),\
v\in Y$, $t\in I$, satisfies the following:

\begin{itemize}
\item for any $v\in Y$ and $t\in I$ operator $-A(t,v)$ is quasi-m-accretive,
that is, it generates a $C_{0}$-semigroup in $X$ and there exists $\beta
=\beta (R)\in \mathbb{R}$ such that%
\begin{equation}
\left\Vert e^{-sA(t,v)}\right\Vert _{X}\leq e^{\beta s},\ s\geq 0,\
\left\Vert v\right\Vert _{Y}\leq R;  \label{A-est}
\end{equation}

\item there exists an isomorphism $Q:Y\rightarrow X$ and $B(t,v)\in \mathcal{%
L(}X,X)$ such that, for all $v\in Y$ and $t\in I$, we have%
\begin{equation}
QA(t,v)Q^{-1}=A(t,v)+B(t,v);  \label{Kato-B}
\end{equation}
the map $I\ni t\mapsto B(t,v)\in X$ is strongly measurable and
\begin{equation}
\lambda =\lambda (R):=\sup_{t\in I}\sup_{v:\left\Vert v\right\Vert _{Y}\leq
R}\left\Vert B(t,v)\right\Vert <\infty ;  \label{B-est}
\end{equation}

\item for any $v\in Y$ and $t\in I$ we have $A(t,v)\in \mathcal{L}(Y,X)$ and
the map
\begin{equation}
I\ni t\mapsto A(t,v)\in \mathcal{L}(Y,X)  \label{A-int}
\end{equation}%
is continuous;

\item there exists $\mu _{A}=\mu _{A}(R)$ such that for all $u\in Y$ and $%
\left\Vert v_{1}\right\Vert _{Y},\left\Vert v_{2}\right\Vert _{Y}\leq R$ we
have%
\begin{equation}
\left\Vert \left( A(t,v_{1})-A(t,v_{2})\right) u\right\Vert _{X}\leq \mu
_{A}\left\Vert v_{1}-v_{2}\right\Vert _{X}\left\Vert u\right\Vert _{Y}.\
\label{Kato-A-est}
\end{equation}
\end{itemize}
\end{condition}

\begin{theorem}
\label{Kato}Let Condition \ref{Kato-cond} hold on the time interval $I=\left[
0,T\right] $. Then for every $v_{0}\in Y$ there exists $T^{\prime
}=T^{\prime }(v_{0})\leq T$ and a unique solution $v\in C\left( [0,T^{\prime
}],Y\right) \cap C^{1}\left( [0,T^{\prime }],X\right) $ of equation (\ref%
{ODE_semi}).
\end{theorem}

\noindent \textbf{Proof.} See \cite[Theorem 6, page 36]{Kato}
\hfill%
$\square $

\begin{remark}
$T^{\prime }$ is an arbitrary number satisfying the following bounds:%
\begin{eqnarray*}
\exp \left( \left( \beta +\lambda \right) T^{\prime }\right) &<&R\left\Vert
v_{0}\right\Vert _{Y}^{-1}, \\
T^{\prime }\exp \left( \beta T^{\prime }\right) &<&R^{-1}\mu _{A}^{-1},
\end{eqnarray*}%
Here the constants $\beta =\beta (R)$, $\lambda =\lambda (R)$ and $\mu
_{A}=\mu _{A}(R)$ are defined in (\ref{A-est}), (\ref{B-est}) and (\ref%
{Kato-A-est}), respectively, see \cite[p. 45]{Kato}. The corresponding
solution of equation (\ref{ODE_semi}) will satisfy the bound $\left\Vert
v(t)\right\Vert _{Y}\leq R$.
\end{remark}

\subsection{Kato's condition for first order differential operators}

We set $X=H^{1},$ $Y=H^{2}$ and $Q=(1-\partial _{x}^{2})^{1/2}$. It is clear
that $Q:H^{2}\rightarrow H^{1}$ is an isometry. We first consider the family
of first order differential operators%
\begin{equation}
\mathcal{A(}y\mathcal{)}=a(y)\partial _{x}+b(y),\ y\in H^{1},  \label{A-op}
\end{equation}%
defined on $H^{2}$, with coefficients $\ a(y)\in H^{3},b(y)\in H^{2},$ $y\in
H^{1}$. We assume that the maps%
\begin{equation}
a:H^{1}\rightarrow H^{3},\ b:H^{1}\rightarrow H^{2}\text{ are Lipschitz
continuous}  \label{ab-lip}
\end{equation}%
and bounded (uniformly in $y$), that is,
\begin{equation}
\sup_{y\in H^{1}}\left\Vert a(y)\right\Vert _{H^{3}}<\infty ,\ \sup_{y\in
H^{1}}\left\Vert b(y)\right\Vert _{H^{2}}<\infty .  \label{coef-bound}
\end{equation}%
It is clear that $\mathcal{A(}y\mathcal{)\in L}(Y,X)$ with the uniformly (in
$y\in H^{1}$) bounded norm.

According to the results of Section \ref{sog} (Lemma \ref{semi-lemma}
below), for any $y\in H^{1}$, there exists a one-parameter $C_{0}$-group in $%
H^{1}$ such that its generator contains $H^{2}$ in its domain and coincides
with $\mathcal{A(}y\mathcal{)}$ on $H^{2}$. We will preserve the notation $%
\mathcal{A(}y\mathcal{)}$ for this operator. Observe that, again by Lemma %
\ref{semi-lemma}, there exists an operator $\mathcal{A}^{(0)}\mathcal{(}y%
\mathcal{)}$ in $L^{2}$, which coincides with $\mathcal{A(}y\mathcal{)}$ on $%
H^{2}$ and generates a one-parameter $C_{0}$-group in $L^{2}$.

\begin{theorem}
\label{Kato-B-lemma}The operator family (\ref{A-op}) satisfies Condition \ref%
{Kato-cond} on the time interval $I=\left[ -\tau ,\tau \right] $, with
arbitrary $R>0$ (appearing in Condition \ref{Kato-cond}) and $\tau >0$.
\end{theorem}

\noindent \textbf{Proof. }(i)\textbf{\ }The first\textbf{\ }part of
Condition \ref{Kato-cond} immediately follows from the results of Section %
\ref{sog} below. Indeed, the fact that $-\mathcal{A(}y\mathcal{)}$ is the
generator of the $C_{0}$-semigroup in $X$ and estimate (\ref{A-est}) follow
from Lemma \ref{semi-lemma} below and the bound (\ref{coef-bound}).

(ii) Condition (\ref{Kato-B}) is essentially proved in \cite[Remark 2.6 b)]%
{CE} for $a(y)=Q^{-2}v$ and $b(y)=2\left( \partial _{x}Q^{-2}v\right) ,\
v=Qy\in L^{2}$, cf. (\ref{A-CH}). The proof does not use the explicit form
of the coefficients. Here we give its main steps adapted to our setting.

We fix $y\in H^{1}$ and use the shorthand notation $\mathcal{A}:=\mathcal{A(}%
y\mathcal{)}$ and $\mathcal{A}^{(0)}:=\mathcal{A}^{(0)}(y)$. The first step
is to prove equality (\ref{Kato-B}) for the operator $\mathcal{A}^{(0)}$ in
the pair of spaces $H^{1}\subset L^{2}$. Denote by $M_{a}$ and $M_{b}$ the
operators of multiplication by $a:=a(y)\in H^{3}$ and $b:=b(y)\in H^{2}$,
respectively. Define an operator $B$ by the equality $Bf:=Q\mathcal{A}%
^{(0)}Q^{-1}f-\mathcal{A}^{(0)}f\mathcal{\ }$for $f\in \mathcal{S}%
:=C^{\infty }\cap L^{2}$. Then on $\mathcal{S}$ we have the equality%
\begin{equation*}
B\mathcal{=}\left[ Q,M_{a}\right] \partial _{x}Q^{-1}+QM_{b}Q^{-1}-M_{b},
\end{equation*}%
because $\partial _{x}Q^{-1}f=Q^{-1}\partial _{x}f$ for $f\in \mathcal{S}$.
The operators $M_{b},\ QM_{b}Q^{-1}$ and $\partial _{x}Q^{-1}$ are bounded
in both spaces $L^{2}$ and $H^{1}$, with
\begin{equation}
\left\Vert QM_{b}Q^{-1}\right\Vert _{\mathcal{L(}L^{2})}=\left\Vert
M_{b}\right\Vert _{\mathcal{L(}H^{1})}\leq \left\Vert b\right\Vert _{H^{1}}
\label{QMQ1}
\end{equation}%
and
\begin{equation}
\left\Vert QM_{b}Q^{-1}\right\Vert _{\mathcal{L(}H^{1})}=\left\Vert
M_{b}\right\Vert _{\mathcal{L(}H^{2})}\leq \left\Vert b\right\Vert _{H^{2}}
\label{QMQ}
\end{equation}%
(because $H^{1}$ and $H^{2}$ are Banach algebras).

It is proved in \cite[Section VII.3.5]{Stein} that the commutator $\left[
Q,M_{a}\right] $ is bounded in $L^{2}$ and%
\begin{equation}
\left\Vert \left[ Q,M_{a}\right] \right\Vert _{\mathcal{L(}L^{2})}\leq
K\left\Vert \partial _{x}a\right\Vert _{H^{1}},  \label{commut-est-L2}
\end{equation}%
for some constant $K>0$. This bound together with (\ref{QMQ1}) implies that $%
B$ is a bounded operator in $L^{2}$ and
\begin{equation*}
\left\Vert B\right\Vert _{\mathcal{L(}L^{2})}\leq K\max \left( \left\Vert
\partial _{x}a\right\Vert _{H^{1}},\left\Vert b\right\Vert _{H^{1}}\right) .
\end{equation*}%
It is proved in \cite[Proposition 2.3 a)]{CE} that $\mathcal{S}$ is a core
for $\mathcal{A}^{(0)}$, which is sufficient for the equality%
\begin{equation}
Q\mathcal{A}^{(0)}Q^{-1}=\mathcal{A}^{(0)}+B  \label{A0}
\end{equation}%
to hold (\cite[Remark 7.1.3.]{Kato}).

We observe that the operator $\mathcal{A}$ coincides with the part of $%
\mathcal{A}^{(0)}$ in $H^{1}$ (\cite[Theorem 4.5.5 and Lemma 5.4.4 ]{Pazy}).
Thus, equality (\ref{Kato-B}) for $\mathcal{A}$ will follow from (\ref{A0})
provided $B\in \mathcal{L}(H^{1})$. As in \cite[Proposition 2.3 a)]{CE}, we
can write%
\begin{eqnarray*}
\left\Vert \left[ Q,M_{a}\right] \right\Vert _{\mathcal{L(}H^{1})}^{2}
&=&\left\Vert \left[ Q,M_{a}\right] Q^{-1}\right\Vert _{\mathcal{L(}%
L^{2},H^{1})}^{2} \\
&\leq &\left\Vert \left[ Q,M_{a}\right] Q^{-1}\right\Vert _{\mathcal{L(}%
L^{2})}^{2}+\left\Vert \partial _{x}\left[ Q,M_{a}\right] Q^{-1}\right\Vert
_{\mathcal{L(}L^{2})}^{2}.
\end{eqnarray*}%
The first term is bounded by $K\left\Vert Q^{-1}\right\Vert _{\mathcal{L(}%
L^{2})}^{2}\left\Vert \partial _{x}a\right\Vert _{H^{1}}^{2}\leq c\left\Vert
\partial _{x}a\right\Vert _{H^{2}}^{2}$, cf. (\ref{commut-est-L2}). For the
second term we have%
\begin{equation*}
\partial _{x}\left[ Q,M_{a}\right] Q^{-1}=QM(\partial
_{x}a)Q^{-1}+M(\partial _{x}a)+\left[ Q,M_{a}\right] \partial _{x}Q^{-1},
\end{equation*}%
which, together with (\ref{QMQ}) applied to the operator $QM(\partial
_{x}a)Q^{-1}$ and a new use of (\ref{commut-est-L2}), leads to the bound%
\begin{equation*}
\left\Vert \partial _{x}\left[ Q,M_{a}\right] Q^{-1}\right\Vert _{\mathcal{L(%
}L^{2})}^{2}\leq c\left\Vert \partial _{x}a\right\Vert _{H^{2}}
\end{equation*}%
for some constant $c>0$, and so%
\begin{equation}
\left\Vert B\right\Vert _{\mathcal{L(}H^{1})}\leq K\max \left( \left\Vert
\partial _{x}a\right\Vert _{H^{2}},\left\Vert b\right\Vert _{H^{2}}\right) ,
\label{B-norm-bound}
\end{equation}%
for a generic constant $K>0$.

Finally, estimate (\ref{B-est}) follows now from assumption (\ref{ab-lip}).

(iii) Condition (\ref{A-int}) trivially holds because $A(v)$ is independent
of $t$. Condition (\ref{Kato-A-est}) can be checked directly using (\ref%
{ab-lip}).

\hfill%
$\square $

\begin{remark}
\label{Kato-B-rem}We observe that (\ref{Kato-B}) remains true if the
coefficients $a$ and $b$ in (\ref{A-op}) are $t$-dependent and such that,
for every $y\in H^{1}$, the right-hand side of (\ref{B-norm-bound}) is
bounded uniformly in $t$. For condition (\ref{A-int}) to hold, it is
sufficient that, for every $y\in H^{1}$, the maps $\mathbb{R}\ni t\mapsto
a(t,y)\in H^{2}$ and $\mathbb{R}\ni t\mapsto b(t,y)\in H^{1}$ are continuous.
\end{remark}

\begin{remark}
In \cite[Remark 2.6 b)]{CE}, the authors took a slightly different path.
They proved Condition \ref{Kato-cond} for the pair $X=L^{2}$ and $Y=H^{1}$,
which implies the existence of a solution of the Camassa-Holm equation (\ref%
{QLE}) in $H^{1}$. Then they showed that the solution actually belongs to $%
H^{2}$ provided the initial condition does so.
\end{remark}

\subsection{Proof of the main result.\label{main-proof}}

In this section we will show that Kato's theory can be applied to the
integral equation (\ref{SCH-int}). Recall that%
\begin{equation}
\widehat{A}(t,v)=U_{t}A(U_{t}^{-1}v)U_{t}^{-1},\ v\in H^{1},\ t\in \mathbb{R}%
,  \label{Atv1}
\end{equation}%
cf. (\ref{Atv}). It has been proved in Proposition \ref{ab-pert} that $%
\widehat{A}(t,v)=a(t,v)\partial _{x}+b(t,v)$ with $a(t,v)\in H^{3}$ and $%
b(t,v)\in H^{2}$. As before, we retain the same notation for the generator
of the corresponding one-parameter $C_{0}$-group in $L^{2}$ (see Lemma \ref%
{semi-lemma} below).

\begin{theorem}
\label{Kato-pert}For any $\tau ,R>0$ (with $R$ as in Condition \ref%
{Kato-cond}), the operator family $\widehat{A}(t,v),t\in \left[ -\tau ,\tau %
\right] ,v\in H^{2}$, satisfies Condition \ref{Kato-cond} with $X=H^{1}$ and
$Y=H^{2}$.
\end{theorem}

\noindent \textbf{Proof. } It is clear that the coefficients $a(t,v)\in
H^{3} $ and $b(t,v)\in H^{2}$ are bounded uniformly in $t$ so that (\ref%
{B-norm-bound}) is satisfied. Also, the Lipschitz condition (\ref{ab-lip})
holds because of the explicit form (\ref{at}), (\ref{bt}) of the
coefficients and uniform in $t\in \left[ -\tau ,\tau \right] $ boundedness
of the group $U_{t}$ in both $\mathcal{L}(X,X)$ and $\mathcal{L}(Y,Y)$ (cf. (%
\ref{U-est0})). Thus, according to Theorem \ref{Kato-B-lemma} and Remark \ref%
{Kato-B-rem}, the operator family $\widehat{A}(t,v),\ t\in \left[ -\tau
,\tau \right] $, satisfies the first two parts of Condition \ref{Kato-cond}
with arbitrary $R$.\textbf{\ }

The continuity condition (\ref{A-int}) is proved in Corollary \ref{a-cont}.
Estimate (\ref{Kato-A-est}) immediately follow from (\ref{Atv1}) as well as
the (uniform in $t\in \left[ -\tau ,\tau \right] $) boundedness of operators
$U_{t}$ in both $X$ and $Y$, cf. (\ref{U-est0}).
\hfill%
$\square $

\begin{remark}[Change of time]
\label{time-change}Let $f:\left[ 0,T\right] \rightarrow \left[ -\tau ,\tau %
\right] $ be a continuous function. It is clear that operator family $%
A_{f}(t,v):=\widehat{A}(f(t),v),\ t\in \left[ 0,T\right] $, satisfies
Condition \ref{Kato-cond}. Moreover, since $\sup_{t\in \left[ 0,T\right]
}\left\Vert U_{f(t)}\right\Vert \leq \sup_{t\in \left[ -\tau ,\tau \right]
}\left\Vert U_{t}\right\Vert $, the constants $\beta ,\lambda $ and $\mu
_{A} $ remain unchanged.
\end{remark}

We return now to the stochastic Camassa-Holm equation (\ref{SPDE}), defined
on the filtered probability space $\left( \Omega ,\mathcal{F},\left(
\mathcal{F}_{t}\right) _{t\geq 0},\mathbb{P}\right) $.

\begin{theorem}
\label{AWK}For any $R>0$ and $z_{0}\in \mathbb{R}$ and each continuous
Brownian path $w(t)$ there exists $\theta >0$ and a unique solution $z\in C(%
\left[ 0,\theta \right] ,H^{2})$, of the integral equation (\ref{SCH-int}),
such that $z(0)=z_{0}$ and $\left\Vert z(t)\right\Vert _{H^{2}}\leq R$, $%
t\in \left[ 0,\theta \right] $.
\end{theorem}

\noindent \textbf{Proof.} Fix $R>0$ and a continuous Brownian path $w(t)$.
Fix in addition $T>0$ and define $\tau =\tau (w):=\inf \left\{
t>0:\left\vert w(t)\right\vert \geq T\right\} $. According to Theorem \ref%
{Kato-pert} and Remark \ref{time-change}, the operator family $\widehat{A}%
(w(t),v),$ $t\in \left[ 0,\tau \right] ,$ satisfies Condition \ref{Kato-cond}
with the constants $\beta ,\lambda $ and $\mu _{A}$ (depending on $R$ and $T$%
).

Next, we choose any $T^{\prime }>0$ such that%
\begin{eqnarray*}
\exp \left( \left( \beta +\lambda \right) T^{\prime }\right) &\leq
&R\left\Vert v_{0}\right\Vert _{H^{2}}^{-1}, \\
T^{\prime }\exp \left( \beta T^{\prime }\right) &<&R^{-1}\mu _{A}^{-1},
\end{eqnarray*}%
and define $\theta :=\min \left\{ \tau ,T^{\prime }\right\} $. Then, by
Theorem \ref{Kato}, there exists a solution $z\in C(\left[ 0,\theta \right]
,H^{2})$ of the integral equation (\ref{SCH-int}), such that $\left\Vert
z(t)\right\Vert _{H^{2}}\leq R$, $t\in \left[ 0,\theta \right] $.
\hfill%
$\square $

\begin{remark}
It is clear that, for any $R>0$, both $\tau $ and $\theta $ are stopping
times.
\end{remark}

\noindent \textbf{Proof of Theorem \ref{theor-main}. }The process\textbf{\ }$%
z(t)$ constructed in Theorem \ref{AWK} satisfies the conditions of Theorem %
\ref{ode-theor} with $Y=H^{2}$ and $\mathfrak{X}=L^{2}$, which implies that $%
y(t):=U_{w(t)}^{-1}z(t),\ t\in \left[ 0,\theta \right] $, is the unique
strong solution of equation (\ref{SPDE}).
\hfill%
$\square $

\section{Auxiliary results\label{aux-res}}

In this section we present some general results used in the main part of the
paper.

\subsection{One-parameter groups generated by first order differential
operators\label{sog}}

The aim of this section is to discuss properties of one-parameter groups in
Sobolev spaces $H^{n}$, $n=0,1,2,...,$ generated by first order differential
operators. We will use the convention $H^{0}=L^{2}$.

We need some preparations. Let $C_{b}^{n}$ be the Banach space of $n$-times
continuously differentiable functions $f:\mathbb{R}\rightarrow \mathbb{R}$
with the norm%
\begin{equation*}
\left\Vert f\right\Vert ^{(n)}:=\max_{m=0,...,n}\sup_{x\in \mathbb{R}%
}\left\vert f^{(m)}(x)\right\vert <\infty ,
\end{equation*}%
where $f^{(m)}$ stands for the $m$-th derivative, $f^{(0)}\equiv f$.

Given a function $g(t,x),\ t,x\in \mathbb{R}$, we will keep the notation $%
g^{(m)}(t,x):=\partial _{x}^{m}g(t,x)$ for the $m$-th derivative w.r.t. $x$.
We will use, where possible, notations $g(t)$ and $g^{(m)}(t)$ for the
mappings $x\mapsto g(t,x)$. and $x\mapsto g^{(m)}(t,x)$, respectively. Thus,
we have $g^{\prime }(t):x\mapsto g^{\prime }(t,x)$.

The following statement is essentially well-known.

\begin{lemma}
\label{lemma-diff}Assume that $\xi \in C_{b}^{n+1}$, $n\geq 0$, and consider
equation%
\begin{equation}
\frac{d}{dt}\psi (t,x)=-\xi (\psi (t,x)),\ \psi (0,x)=x,\ x\in \mathbb{R}.
\label{diff-eq}
\end{equation}%
Then:\newline
(i) there exist a unique solution $\psi (t)$, $t\in \mathbb{R}$, of (\ref%
{diff-eq}); it satisfies the estimate%
\begin{equation}
\left\vert \psi (t,x)\right\vert \leq c_{1}e^{c_{2}\left\vert t\right\vert
}\left( \left\vert x\right\vert +c_{3}\right) ;  \label{phi-est}
\end{equation}%
moreover, $\phi (t):=\psi (t)-\mathrm{id}\in L^{2}$.\newline
(ii) the solution $\psi $ is $x$differentiable; moreover, for any $t\in
\mathbb{R}$, the derivative $\psi ^{(1)}(t)\in C_{b}^{n}$, and the following
estimate holds:%
\begin{equation}
\left\Vert \psi ^{(1)}(t)\right\Vert ^{(n)}\leq e^{c_{4}\left\vert
t\right\vert }.  \label{est-diff}
\end{equation}%
The map
\begin{equation*}
\mathbb{R}_{+}\ni t\mapsto \psi ^{(1)}(t)\in C_{b}^{n}
\end{equation*}%
is continuously differentiable.\newline
Here $c_{1},c_{2},c_{3},c_{4}>0$ are some constants depending only on $%
\left\Vert \xi \right\Vert ^{(n)}$.
\end{lemma}

\noindent \textbf{Proof.}

(i) For any fixed $x\in \mathbb{R}$, equation (\ref{diff-eq}) has a solution
because its right-hand side is globally Lipschitz. Estimate (\ref{phi-est})
follows in a standard way from the Gronwall inequality. Since $\phi (t)$
satisfies equation%
\begin{equation*}
\frac{d}{dt}\phi (t,x)=-\widetilde{\xi }(\phi (t,x)),\ \phi (0,x)=0,
\end{equation*}%
where $\widetilde{\xi }(x):=\xi (x)+x$ and is globally Lipschitz in $L^{2}$,
the result follows.

(ii) Consider the linear operator $\widehat{\xi }(t)$ acting on functions $u:%
\mathbb{R\rightarrow R}$ by multiplication by $\xi ^{(1)}(\psi (t,\cdot ))$,
that is,
\begin{equation*}
(\widehat{\xi }(t)u)(x):=\xi ^{(1)}(\psi (t,x))u(x).
\end{equation*}%
A direct calculation shows that $\widehat{\xi }$ is a bounded operator in $%
C_{b}^{n}$ with norm
\begin{equation*}
\left\Vert \widehat{\xi }\right\Vert _{\mathcal{L}(C_{b}^{n-1})}=a_{n}\left%
\Vert \xi \right\Vert ^{(n)},\ 0<a_{n}<\infty .
\end{equation*}%
It is immediate that $\psi ^{(1)}(t)$ solves the equation%
\begin{equation*}
\frac{d}{dt}\psi ^{(1)}(t)=-\widehat{\xi }(t)\psi ^{(1)}(t),\ \psi
^{(1)}(0)=1.
\end{equation*}%
This equation has a unique solution in $C_{b}^{n}$, which satisfies (\ref%
{est-diff}) and is continuously differentiable in $t$.
\hfill%
$\square $

\begin{remark}
In particular, Lemma \ref{lemma-diff} implies in a standard way that $\left(
\psi (t)\right) _{t\in \mathbb{R}}$ is a one-parameter group of $C^{n+1}$%
-diffeomorphisms of $\mathbb{R}^{1}$, generated by the vector field $-\xi
\partial _{x}$.
\end{remark}

Let us introduce an operator family $U_{t}^{\xi }$, $t\in \mathbb{R}$, by
the formula $U_{t}^{\xi }f=f(\psi (t)),\ f:\mathbb{R}\rightarrow \mathbb{R}$.

\begin{lemma}
\label{lemma-U}Assume that $\xi \in C_{b}^{n+1}$, $n\geq 0$. Then $%
U_{t}^{\xi },\ t\in \mathbb{R}$, is a strongly continuous one-parameter
group of bounded operators in $H^{n}$ such that%
\begin{equation}
\left\Vert U_{t}^{\xi }\right\Vert _{\mathcal{L}(H^{n})}\leq
c_{1}e^{c_{2}\left\vert t\right\vert },\ t\in \mathbb{R},  \label{u-est}
\end{equation}%
for some positive constants $c_{1},c_{2}<\infty $ (depending only on $n$ and
$\left\Vert \xi \right\Vert ^{(n+1)}$). In the case of $n=0$ and $n=1$ we
can take $c_{1}=1$. The domain of the generator $D_{0}$ of $U_{t}^{\xi }$
contains $H^{n+1}$ and one has $D_{0}=\xi \partial _{x}$ on $H^{n+1}$.
\end{lemma}

\noindent \textbf{Proof.} In this proof, $c,c_{1},c_{2},...$ will stand for
universal positive constants (depending only on $n$ and $\left\Vert \xi
\right\Vert ^{(n+1)}$).

1) Let us prove that the operators $U_{t}^{\xi },\ t\in \mathbb{R},$ are
bounded in $H^{n}$.

Consider first the case of $n=0$. Then, for $f\in H^{0}\equiv L^{2}$, we have%
\begin{equation*}
\left\Vert U_{t}^{\xi }f\right\Vert _{L^{2}}^{2}=\int f(\psi
(t,x))^{2}dx=\sup_{x\in \mathbb{R}}\left\vert \psi ^{(1)}(-t,x)\right\vert
^{2}\left\Vert f\right\Vert _{L^{2}}^{2}\leq e^{2c_{4}\left\vert
t\right\vert }\left\Vert f\right\Vert _{L^{2}}^{2},
\end{equation*}%
cf. (\ref{est-diff}), and estimate (\ref{u-est}) holds with $c_{1}=1$.

Let now $n\geq 1$. By Fa\`{a} di Bruno's theorem for any $k=1,2,...n$ we
have
\begin{equation}
\partial _{x}^{k}f(\psi (t,x))=\sum_{m=1}^{k}f^{(m)}(\psi (t,x))B_{k.m}(\psi
^{(1)}(t,x),...,\psi ^{(k-m+1)}(t,x)),  \label{FdB}
\end{equation}%
where $B_{k.m}$ is the exponential Bell polynomial. Thus, we have%
\begin{multline}
\left\vert \partial _{x}^{k}f(\psi (t,x))\right\vert ^{2}  \label{FdB-est} \\
\leq k\max_{m=1,...,k}\sup_{x\in \mathbb{R}}\left\vert B_{k.m}(\psi
^{(1)}(t,x),...,\psi ^{(k-m+1)}(t,x))\right\vert
^{2}\sum_{m=1}^{k}\left\vert f^{(m)}(\psi (t,x))\right\vert ^{2}.
\end{multline}

It follows from Lemma \ref{lemma-diff} (ii) that, for any $m=1,...,n$ we
have $\psi ^{(m)}(t)\in C_{b}$ and $\sup_{x}\left\vert \psi
^{(m)}(t,x)\right\vert \leq e^{tc_{4}}$. Thus we obtain the estimate%
\begin{equation}
\left\Vert f(\psi (t))\right\Vert _{H^{n}}\leq c_{1}e^{c_{2}\left\vert
t\right\vert }\left\Vert f\right\Vert _{H^{n}},  \label{H-bound}
\end{equation}%
for some constants $0<c_{1},c_{2}<\infty $ (depending only on $n$ and $%
\left\Vert \xi \right\Vert ^{(n)}$), which implies (\ref{u-est}).

By observe that in the case of $n=1$ formula (\ref{FdB-est}) gets the form%
\begin{equation*}
\left\vert \partial _{x}f(\psi (t,x))\right\vert ^{2}\leq \sup_{x\in \mathbb{%
R}}\left\vert \psi ^{(1)}(t)(x)\right\vert ^{2}\left\vert f^{(1)}(\psi
(t,x))\right\vert ^{2},
\end{equation*}%
it follows from (\ref{est-diff}) that%
\begin{multline*}
\left\Vert f(\psi (t))\right\Vert _{H^{1}}^{2}\leq \int \left[ \left\vert
f(\psi (t,x))\right\vert ^{2}+\sup_{x\in \mathbb{R}}\left\vert \psi
^{(1)}(t,x)\right\vert ^{2}\left\vert f^{(1)}(\psi (t,x))\right\vert ^{2}%
\right] dx \\
\leq \sup_{x\in \mathbb{R}}\left\vert \psi ^{(1)}(t,x)\right\vert
^{2}\sup_{x\in \mathbb{R}}\left\vert \psi ^{(1)}(-t,x)\right\vert
^{2}\left\Vert f\right\Vert _{H^{1}}^{2} \\
\leq e^{4\left\vert t\right\vert c_{4}}\left\Vert f\right\Vert _{H^{1}}^{2},
\end{multline*}%
and estimate (\ref{H-bound}) holds with $c_{1}=1$ and $c_{2}=4c_{4}$.

2) The fact that the operators $U_{t}^{\xi }\in \mathcal{L}(H^{n}),\ t\in
\mathbb{R},$ form a one-parameter group, that is, $U_{t_{1}}^{\xi
}U_{t_{2}}^{\xi }=U_{t_{1}+t_{2}}^{\xi }$, $t_{1},t_{2}\in \mathbb{R}$, and $%
U_{0}^{\xi }=I$, follows in a standard way from the group properties of the
family of diffeomorphisms $\psi (t)$, $t\in \mathbb{R}$.

3) Now we shall prove that the map $\mathbb{R}\ni t\mapsto U_{t}^{\xi }$ is
strongly continuous. Let $f\in C_{0}^{\infty }$. Observe that $\psi
^{(m)}(t,x)\rightarrow x^{(m)}=\left\{
\begin{array}{c}
1,\ m=1 \\
0,\ m\geq 2%
\end{array}%
\right. ,\ t\rightarrow 0$, uniformly on compact sets. Thus for the r.h.s.
of (\ref{FdB}) we have%
\begin{multline*}
\sum_{m=1}^{k}f^{(m)}(\psi (t,x))B_{k.m}(\psi ^{(1)}(t,x),...,\psi
^{(k-m+1)}(t,x)) \\
\rightrightarrows
\sum_{m=1}^{k}f^{(m)}(x))B_{k.m}(x^{(1)},...,x^{(k-m+1})=f^{(k)}(x),\
t\rightarrow 0,
\end{multline*}%
where $\rightrightarrows $ stands for the uniform convergence in $x\in
\mathbb{R}$. The last equality holds because $x^{(m)}=\left\{
\begin{array}{c}
1,\ m=1 \\
0,\ m\geq 2%
\end{array}%
\right. $ and thus $B_{k.m}(x^{(1)},...,x^{(k-m+1})=\left\{
\begin{array}{c}
1,\ m=k \\
0,\ m\leq k-1%
\end{array}%
\right. $. Therefore $\partial _{x}^{k}f(\psi (t))\overset{L_{2}}{%
\rightarrow }f^{(k)},\ t\rightarrow 0$, for any $k\leq n$, which implies the
convergence $f(\psi (t))\overset{H^{n}}{\rightarrow }f,\ t\rightarrow 0$.

Let now $u\in H^{n}$. We have the estimate%
\begin{multline*}
\left\Vert U_{t}^{\xi }u-u\right\Vert _{H^{k}}\leq \left\Vert U_{t}^{\xi
}u-U_{t}^{\xi }f\right\Vert _{H^{k}}+\left\Vert U_{t}^{\xi }f-f\right\Vert
_{H^{k}}+\left\Vert f-u\right\Vert _{H^{k}} \\
\leq \left\Vert U_{t}^{\xi }f-f\right\Vert _{H^{k}}+c\left\Vert
f-u\right\Vert _{H^{k}},
\end{multline*}%
and the required result follows from the fact that $C_{0}^{\infty }$ is
dense in $H^{n}$.

3) Let us prove that $t\mapsto U_{t}^{\xi }u\in H^{n}$ is differentiable for
$u\in H^{n+1}$.

Let $v\in C_{0}^{\infty }$. Formula (\ref{FdB}) implies that $t\mapsto
\partial _{x}^{k}v(\psi (t,x))$ is continuously differentiable for any $x\in
\mathbb{R}$. Denote
\begin{equation*}
F(x):=\frac{d}{dt}\partial _{x}^{k}v(\psi (t,x))_{t=0}.
\end{equation*}%
Then%
\begin{equation*}
\frac{\partial _{x}^{k}v(\psi (t,x))-v^{(k)}(x)}{t}=\partial _{x}^{k}\frac{%
v(\psi (t,x))-v(x)}{t}\rightrightarrows F(x),\ t\rightarrow 0,
\end{equation*}%
because $f$ has compact support, and so%
\begin{equation}
\frac{v(\psi (t))-v}{t}\overset{H^{n}}{\longrightarrow }\frac{d}{dt}v(\psi
(t))_{t=0},\ t\rightarrow 0.  \label{H-der}
\end{equation}

We will prove now that (\ref{H-der}) holds for any $u\in H^{n+1}$. Set $%
g(t,x)=\partial _{x}^{k}u(\psi (t,x))$ and $f(t,x)=\partial _{x}^{k}v(\psi
(t,x))$. Then (denoting the derivative w.r.t. the first variable by "dot")
we obtain%
\begin{multline}
\dot{g}(s,x)=\frac{d}{ds}\partial _{x}^{k}u(\psi (s,x))=\partial _{x}^{k}%
\frac{d}{ds}u(\psi (s,x)) \\
=\partial _{x}^{k}\left[ \xi (x)\partial _{x}u(\psi (s,x))\right]
=\sum_{m=0}^{k}\left(
\begin{array}{c}
k \\
m%
\end{array}%
\right) \xi ^{(m)}(x)\partial _{x}^{k-m}\partial _{x}u(\psi (s,x)).
\label{g-dot}
\end{multline}%
In particular,
\begin{equation}
\dot{g}(0,x)=\sum_{m=0}^{k}\left(
\begin{array}{c}
k \\
m%
\end{array}%
\right) \xi ^{(m)}(x)\partial _{x}^{k-m+1}u(x).  \label{g-dot-0}
\end{equation}
Of course, similar formulae hold for $f$.

Thus, applying Cauchy inequality, we obtain%
\begin{multline*}
\int_{\mathbb{R}}\left\vert \frac{1}{t}\int_{0}^{t}\left( \dot{g}(s,x)-\dot{f%
}(s,x)\right) ds\right\vert ^{2}dx \\
\leq \frac{1}{t}\int_{\mathbb{R}}\int_{0}^{t}\left\vert \dot{g}(s,x)-\dot{f}%
(s,x)\right\vert ^{2}dsdx \\
\leq \sum_{m}\left(
\begin{array}{c}
k \\
m%
\end{array}%
\right) \frac{1}{t}\int_{0}^{t}\int_{\mathbb{R}}\left\vert \xi
^{(m)}(x)\right\vert \left\vert \partial _{x}^{k-m+1}\left( u(\psi
(s,x))-v(\psi (s,x))\right) \right\vert ^{2}dxds \\
\leq c\left\Vert \xi \right\Vert _{C_{b}^{k}}^{2} \\
\sum_{m}\left(
\begin{array}{c}
k \\
m%
\end{array}%
\right) \frac{1}{t}\int_{0}^{t}\int_{\mathbb{R}}\left\vert u^{(k-m+1)}(\psi
(s,x))-v^{(k-m+1)}(\psi (s,x))\right\vert ^{2}dxds.
\end{multline*}%
The last inequality is due to the formulae (\ref{FdB}) and (\ref{H-bound}).
Taking into account that $\int p(\psi (s,x))dx=\int \left\vert \partial
_{x}\psi (s,x)^{-1}\right\vert p(x)dx$ for any integrable function $p$ we
obtain%
\begin{equation*}
\int_{\mathbb{R}}\left\vert \frac{1}{t}\int_{0}^{t}\left( \dot{g}(s,x)-\dot{f%
}(s,x)\right) ds\right\vert ^{2}dx\leq c_{1}\left\Vert u-v\right\Vert
_{H^{k+1}}^{2}.
\end{equation*}%
Observe that (\ref{g-dot-0}) implies that%
\begin{equation*}
\int_{\mathbb{R}}\left\vert \dot{g}(0,x)-\dot{f}(0,x)\right\vert ^{2}dx\leq
c_{2}\left\Vert u-v\right\Vert _{H^{k+1}}^{2}.
\end{equation*}%
The following general relation holds for any $t>0,\ x\in \mathbb{R}$:%
\begin{multline*}
\frac{g(t,x)-g(0,x)}{t}-\dot{g}(0,x)=\frac{1}{t}\int_{0}^{t}\dot{g}(s,x)ds-%
\dot{g}(0,x) \\
=\frac{1}{t}\int_{0}^{t}\left( \dot{g}(s,x)-\dot{f}(s,x)\right) ds+\left[
\frac{1}{t}\int_{0}^{t}\dot{f}(s,x)ds-\dot{f}(0,x)\right] \\
+\left[ \dot{f}(0,x)-\dot{g}(0,x)\right] \\
=\frac{1}{t}\int_{0}^{t}\left( \dot{g}(s,x)-\dot{f}(s,x)\right) ds+\left[
\frac{f(t,x)-f(0,x)}{t}-\dot{f}(0,x)\right] \\
+\left[ \dot{f}(0,x)-\dot{g}(0,x)\right] .
\end{multline*}%
Recalling that $g(t,x)=\partial _{x}^{k}u(\psi (t,x))$ and $f(t,x)=\partial
_{x}^{k}v(\psi (t,x))$, we obtain%
\begin{multline*}
\int_{\mathbb{R}}\left\vert \partial _{x}^{k}\left( \frac{u(\psi (t,x))-u(x)%
}{t}-\frac{d}{dt}u(\psi (t,x))_{t=0}\right) \right\vert ^{2}dx \\
=\int_{\mathbb{R}}\left\vert \frac{g(t,x)-g(0,x)}{t}-\dot{g}(0,x)\right\vert
^{2}dx \\
\leq c_{1}\left\Vert u-v\right\Vert _{H^{k+1}}^{2}+c_{3}\left\Vert \frac{%
v(\psi (t))-v}{t}-\frac{d}{dt}v(\psi (t))_{t=0}\right\Vert _{H^{k}}^{2} \\
+c_{2}\left\Vert u-v\right\Vert _{H^{k+1}}^{2}.
\end{multline*}%
This estimate holds for all $k\leq n$, which implies that%
\begin{multline*}
\left\Vert \frac{u(\psi (t))-u}{t}-\frac{d}{dt}u(\psi (t))_{t=0}\right\Vert
_{H^{n}}^{2}\leq c_{4}\left\Vert u-v\right\Vert _{H^{n+1}}^{2} \\
+c_{3}\left\Vert \frac{v(\psi (t))-v}{t}-\frac{d}{dt}v(\psi
(t))_{t=0}\right\Vert _{H^{n}}^{2},
\end{multline*}%
and the result follows from (\ref{H-der}) and the fact that $C_{0}^{\infty }$
is dense in $H^{n+1}$.

\hfill%
$\square $

We will use the following well-known result.

\begin{theorem}
\label{pazy-th}(\cite[Theorem 3.1.1.]{Pazy}) Let $X$ be a Banach space and
let $A$ be the infinitesimal generator of a $C_{0}$ semigroup $T(t)$ on $X$,
satisfying $\left\Vert T(t)\right\Vert _{X}\leq Me^{\omega t}$ for some
positive constants $M$ and $\omega $. If $B$ is a bounded linear operator on
$X$ then $A+B$ is the infinitesimal generator of a $C_{0}$ semigroup $S(t)$
on $X$, satisfying
\begin{equation}
\left\Vert S(t)\right\Vert _{X}\leq Me^{\left( \omega +M\left\Vert
B\right\Vert _{X}\right) t}.  \label{S-est}
\end{equation}
\end{theorem}

Let us now define for $\eta \in C_{b}^{n}$ the operator $D=D_{0}+\eta $ , $%
Dom(D)=Dom(D_{0})$, so that $D=\xi \partial _{x}+\eta $ on $H^{n+1}$, $n\in
\mathbb{N}$.

\begin{lemma}
\label{semi-lemma}Assume that $\xi \in C_{b}^{n+1}$ and $\eta \in C_{b}^{n}$%
, $n\geq 0$. Then $D$ generates a strongly continuous one-parameter group $%
\left( U_{t}^{\xi ,\eta }\right) _{t\in \mathbb{R}}$ in $H^{n}$, which
satisfies the estimate%
\begin{equation}
\left\Vert U_{t}^{\xi ,\eta }\right\Vert _{\mathcal{L}(H^{n})}\leq
C_{1}e^{C_{2}\left\vert t\right\vert },\ t\in \mathbb{R}  \label{U-est}
\end{equation}%
for some positive constants $C_{1},C_{2}$ (depending only on $n$ and $%
\left\Vert \xi \right\Vert ^{(n)},\left\Vert \eta \right\Vert ^{(n)}$). In
the case where $n=0$ or $n=1$ we can take $C_{1}=1$.
\end{lemma}

\noindent \textbf{Proof.} We observe that the operator $D-D_{0}=\eta $ is
bounded in $H^{n}$. The statement follows now from Theorem \ref{pazy-th}.
\hfill%
$\square $

The group $U_{t}^{\xi ,\eta }$ has the following explicit form.

\begin{lemma}
\label{U}For any $f\in H^{n}$ we have%
\begin{equation}
U_{t}^{\xi ,\eta }f(x)=e^{c(t,x)}f(\varphi _{-t}(x)),\ t,x\in \mathbb{R},
\label{U-expl}
\end{equation}%
where $\left( \varphi _{t}\right) _{t\in \mathbb{R}}$ is the diffeomorphism
group generated by the vector field $\xi \partial _{x}$ and
\begin{equation}
c(t,x)=\int_{0}^{t}\eta (\varphi _{s-t}(x))ds,\ t,x\in \mathbb{R}
\label{c-expl}
\end{equation}
\end{lemma}

\noindent \textbf{Proof.} A direct calculation show that the function $%
u(t,x):=e^{c(t,x)}f(\varphi _{-t}(x))$, $t,x\in \mathbb{R}$, is a solution
of the initial value problem $u_{t}=Du,\ u(0,x)=f(x)$, if and only if $%
c(t,x) $ satisfies%
\begin{equation*}
c_{t}=\xi c_{x}+\eta ,\ c(0,x)=0.
\end{equation*}%
Formula (\ref{c-expl}) can be obtained by the method of characteristics or
checked directly (as in fact formula (\ref{U-expl}) itself).
\hfill%
$\square $

\subsection{From SDE to ODE \label{S-O}}

Let us consider a pair of densely embedded Hilbert spaces $Y\subset
\mathfrak{X}$, a continuous map $F:Y\rightarrow \mathfrak{X}$ and a linear
(unbounded) operator $D$ in $\mathfrak{X}$ such that $Y\subset dom(D^{2})$.
Assume that $T>0$ is fixed. Our aim is to study the stochastic differential
equation%
\begin{equation}
dy(t)+F(y(t))dt+Dy(t)\circ dw(t)=0,\ t\in \left[ 0,T\right] ,
\label{SDE-gen}
\end{equation}%
where $w$ is an $\mathbb{R}$-valued Wiener process on a filtered probability
space \newline
$\left( \Omega ,\mathcal{F},\left( \mathcal{F}_{t}\right) _{t\geq 0},\mathbb{%
P}\right) $, and $\circ $ means the Stratonovich stochastic differential. We
suppose without loss of generality that all trajectories of $w$ are
continuous.

\begin{definition}
\label{strong-sol-abst}A strong solution of equation (\ref{SDE-gen}) is a $Y$%
-valued continuous process $y(t)$, $t\in \left[ 0,\theta \right] $, where $%
\theta $ is a stopping time, $0<\theta \leq T$, such that the equality%
\begin{equation*}
y(t\wedge \theta )=y_{0}+\int_{0}^{t\wedge \theta }F(y(s))ds+\frac{1}{2}%
\int_{0}^{t\wedge \theta }D^{2}y(s)ds+\int_{0}^{t\wedge \theta }Dy(s)dw(s),\
\end{equation*}%
$y_{0}\in Y,\ t\geq 0,$ is satisfied in $\mathfrak{X}$, $\mathbb{P}$-a.s.
\end{definition}

Assume now that $D$ is the generator of a one-parameter $C_{0}$ group $%
\left\{ U(t)\right\} _{t\in \mathbb{R}}$ in $\mathfrak{X}$, which leaves $Y$
invariant and satisfies the estimates%
\begin{equation}
\left\Vert U(t)\right\Vert _{\mathcal{L(}\mathfrak{X})}\leq Me^{m\left\vert
t\right\vert },\ \left\Vert U(t)\right\Vert _{\mathcal{L(}Y)}\leq
Me^{m\left\vert t\right\vert }  \label{U-esttt}
\end{equation}%
for some positive constants $M$ and $m$. Let us define a (random) map $%
\widehat{F}:\mathbb{R}_{+}\times Y\rightarrow \mathfrak{X}$
\begin{equation*}
\widehat{F}(t,z):=U(w(t))F(U^{-1}(w(t))z),\ z\in Y,\ t\geq 0.
\end{equation*}%
Obviously, for all $t\geq 0$, $\widehat{F}(t,\cdot )$ is a continuous map $%
Y\rightarrow \mathfrak{X}$. Observe also that the map $\mathbb{R\ni }%
t\mapsto \widehat{F}(t,z)\in \mathbb{R}$ is continuous for any trajectory $%
w(t)$ and $z\in Y$. Consider the (random) integral equation

\begin{theorem}
\label{ode-theor}Assume that $\theta $ is a stopping time, $0<\theta \leq T$%
. Let $z(t)$, $t\in \left[ 0,\theta \right] $, be a continuous $Y$-valued
process such that $\mathbb{E}\int_{0}^{\theta }\left\Vert z(s)\right\Vert
_{Y}^{4}ds<\infty $. Then $z(t)$ satisfies the random integral equation%
\begin{equation}
z(t)=z(0)-\int_{0}^{t}\widehat{F}(s,z(s))ds,\ z_{0}\in Y,\ t\in \left[
0,\theta \right] ,  \label{ODE-int}
\end{equation}%
if and only if the process $y(t):=U(-w(t))z(t)\in Y$, $t\in \left[ 0,\theta %
\right] $, is a strong solution of (\ref{SDE-gen}).
\end{theorem}

To prove this theorem, we first need the following general result, which
follows by an application of the It\^{o} formula for Hilbert space valued
functions.

\begin{lemma}
\label{ito-lemma}Assume that $\theta $ is a stopping time, $0<\theta \leq T$%
. Let $\chi (t)$, $t\in \left[ 0,\theta \right] $, be a progressively
measurable $\mathfrak{X}$-valued random process. Define a process $Z(t),\
t\in \left[ 0,\theta \right] $, by the formula%
\begin{equation}
Z(t):=Z_{0}-\int_{0}^{t}\chi (s)ds,\ Z_{0}\in Y,  \label{Ito-z}
\end{equation}%
and assume that $Z(t)\in Y$ for all $t\in \left[ 0,\theta \right] $, and%
\begin{equation}
\mathbb{E}\int_{0}^{\theta }\left\Vert Z(s)\right\Vert _{Y}^{4}ds<\infty .
\label{chi}
\end{equation}
Set%
\begin{equation*}
y(t):=U(-w(t))Z(t)\in Y,\ t\in \left[ 0,\theta \right] .
\end{equation*}%
Then $y(t)$ satisfies the equation%
\begin{multline}
y(t\wedge \theta )=\int_{0}^{t\wedge \theta }U(-w(s))\chi
(s)ds+\int_{0}^{t\wedge \theta }DU(-w(s))Z(s)\circ dw(s) \\
=\int_{0}^{t\wedge \theta }U(-w(s))\chi (s)ds+\int_{0}^{t\wedge \theta
}Dy(s)\circ dw(s),\ t\geq 0,  \label{Ito-y}
\end{multline}%
in $\mathfrak{X}$.
\end{lemma}

\noindent \textbf{Proof\textbf{\ of Lemma \ref{ito-lemma}}. }Consider a map%
\textbf{\ }%
\begin{equation*}
f:K:=\mathbb{R\times }Y\ni (\tau ,y)\mapsto f((\tau ,y))\in \mathfrak{X}
\end{equation*}%
and assume that $f\in C^{2}(K,\mathfrak{X})$ (= the space of twice
continuously differentiable maps $K\rightarrow \mathfrak{X}$). Then $\frac{%
\partial f}{\partial y}((\tau ,y))\in \mathcal{L}(Y,\mathfrak{X})$, $\frac{%
\partial f}{\partial \tau }((\tau ,y))\in \mathfrak{X}$ and $\frac{\partial
^{2}f}{\partial \tau ^{2}}((\tau ,y))\in \mathfrak{X}$. Observe that $\frac{%
\partial f}{\partial \tau }((\tau ,y))$ can be identified with a
(Hilbert-Schmidt) operator $\mathbb{R}\rightarrow \mathfrak{X}$ acting on $%
h\in \mathbb{R}$ by%
\begin{equation*}
\frac{\partial f}{\partial \tau }((\tau ,y))h:=h\frac{\partial f}{\partial
\tau }((\tau ,y))\in \mathfrak{X,}
\end{equation*}%
with the norm equal to $\left\Vert \frac{\partial f}{\partial \tau }((\tau
,y))\right\Vert _{\mathfrak{X}}$.

Define the stochastic process $\xi (t)=\left( w(t),Z(t)\right) $, $t\in %
\left[ 0,\theta \right] $, in $K.$ It is a $K$-valued It\^{o} process such
that%
\begin{equation*}
d\xi (t)=\alpha (t)dt+\beta dw(t),\ t\in \left[ 0,\theta \right] ,
\end{equation*}%
where $\alpha :\left[ 0,\theta \right] \ni t\mapsto \left( 0,\chi (t)\right)
\in K$ and $\beta :\mathbb{R}\rightarrow K$ is a (Hilbert-Schmidt) operator
acting on $h\in \mathbb{R}$ as $\beta h:=h\cdot \left( 1,0_{Y}\right) $ ($%
=(h,0_{Y})$). Here $0_{Y}$ stands for the zero element of $Y$.

Assume in addition that%
\begin{equation}
\mathbb{E}\int_{0}^{\theta }\left\Vert \left( \frac{\partial }{\partial \tau
}f\right) (\xi (s))\right\Vert _{\mathfrak{X}}^{2}ds<\infty \text{ and }%
\mathbb{E}\int_{0}^{\theta }\left\Vert \left( \frac{\partial ^{2}}{\partial
\tau ^{2}}f\right) (\xi (s))\right\Vert _{\mathfrak{X}}^{2}ds<\infty .
\label{ksi-bound}
\end{equation}%
Define now a process $\xi _{\theta }(t)$, $t\geq 0$, by setting%
\begin{equation*}
\xi _{\theta }(t)=\xi (0)+\int_{0}^{t}\alpha (s)\mathbf{1}_{[0,\theta
]}(s)ds+\int_{0}^{t}\beta \mathbf{1}_{[0,\theta ]}(s)dw(s),\ t\geq 0.
\end{equation*}%
It is clear that $\xi _{\theta }(t)=\xi (t)$ for $t\in \lbrack 0,\theta ]$.

It follows then from the general It\^{o} formula in Hilbert spaces, see e.g.
\cite[Theorem VII.1.2]{DF}, that $f(\xi (t))$ is an $\mathfrak{X}$-valued It%
\^{o} process such that%
\begin{multline}
f(\xi _{\theta }(t))=\int_{0}^{t}\frac{\partial f}{\partial y}(\xi _{\theta
}(s))\chi (s)\mathbf{1}_{[0,\theta ]}(s)ds+\frac{1}{2}\int_{0}^{t}\frac{%
\partial ^{2}f}{\partial \tau ^{2}}(\xi _{\theta }(s))ds  \label{Ito-f} \\
+\int_{0}^{t}\frac{\partial f}{\partial \tau }(\xi _{\theta }(s))dw(s),\
t\geq 0.
\end{multline}%
Here $\frac{\partial f}{\partial y}(\xi (s))$ is a bounded operator $%
Y\rightarrow \mathfrak{X}$ and $\frac{\partial f}{\partial \tau }(\xi (s))$
is a (Hilbert-Schmidt) operator $\mathbb{R}\rightarrow \mathfrak{X}$ acting
on $h\in \mathbb{R}$ by%
\begin{equation*}
\frac{\partial f}{\partial \tau }(\xi (s))h:=h\frac{\partial f}{\partial
\tau }(\xi (s))\in \mathfrak{X}.
\end{equation*}%
Finally, $\frac{\partial ^{2}f}{\partial \tau ^{2}}(\xi (s))$ can be
identified with an element of $\mathfrak{X}$.

Set now
\begin{equation}
f((\tau ,y))=U(-\tau )y  \label{f-num}
\end{equation}
so that $y(t)=f(\xi (t))$. Taking into account that $Y\subset Dom(D^{2})$ we
deduce that $f\in C^{2}(K,\mathfrak{X})$ and
\begin{equation*}
\frac{\partial f}{\partial \tau }(\tau ,y)=-DU(-\tau )y,\ \frac{\partial
^{2}f}{\partial \tau ^{2}}(\tau ,y)=D^{2}U(-\tau )y,\ \frac{\partial f}{%
\partial y}(\tau ,y)=U(-\tau ).
\end{equation*}%
It follows now from (\ref{Ito-f}) that%
\begin{multline*}
f(\xi _{\theta }(t))=\int_{0}^{t}U(-w(s))\chi (s)\mathbf{1}_{[0,\theta ]}ds+%
\frac{1}{2}\int_{0}^{t}D^{2}U(-w(s))Z(s)\mathbf{1}_{[0,\theta ]}ds \\
-\int_{0}^{t}DU(-w(s))Z(s)\mathbf{1}_{[0,\theta ]}dw(s),\ t\geq 0,
\end{multline*}%
which implies (\ref{Ito-y}).

Now it is only left to prove (\ref{ksi-bound}), which is equivalent to the
pair of inequalities
\begin{equation}
\mathbb{E}\int_{0}^{\theta }\left\Vert DU(-w(s))Z(s)\right\Vert _{\mathfrak{X%
}}^{2}ds<\infty \text{, }\mathbb{E}\int_{0}^{\theta }\left\Vert
D^{2}U(-w(s))Z(s)\right\Vert _{\mathfrak{X}}^{2}ds<\infty .
\label{ito-bound}
\end{equation}%
Observe that both $D$ and $D^{2}$ are bounded operators from $Y$ to $%
\mathfrak{X}$, so that (\ref{ito-bound}) becomes equivalent to the bound%
\begin{equation}
\mathbb{E}\int_{0}^{\theta }\left\Vert U(-w(s))Z(s)\right\Vert
_{Y}^{2}ds<\infty .  \label{ito-bound1}
\end{equation}%
Then, taking into account that $e^{\lambda \left\vert x\right\vert }\leq
e^{\lambda x}+e^{-\lambda x}$ and $\mathbb{E}e^{\lambda w(s)}=e^{\frac{1}{2}%
\lambda ^{2}s}$, we obtain the bound%
\begin{multline*}
\left( \mathbb{E}\int_{0}^{\theta }\left\Vert U(-w(s))Z(s)\right\Vert
_{Y}^{2}ds\right) ^{2}\leq \mathbb{E}\int_{0}^{\theta }\left\Vert
U(-w(s))\right\Vert ^{4}ds~\mathbb{E}\int_{0}^{\theta }\left\Vert
Z(s)\right\Vert _{Y}^{2}ds \\
\leq M^{4}\int_{0}^{T}\mathbb{E}e^{4m\left\vert w(s)\right\vert }ds~\mathbb{E%
}\int_{0}^{\theta }\left\Vert Z(s)\right\Vert _{Y}^{2}ds \\
\leq 2M^{4}\int_{0}^{T}\left( \mathbb{E}e^{4mw(s)}\right) ds~\mathbb{E}%
\int_{0}^{t}\left\Vert Z(s)\right\Vert _{Y}^{4}ds \\
\leq 2M^{4}\int_{0}^{T}e^{8m^{2}s}ds~\mathbb{E}\int_{0}^{\theta }\left\Vert
Z(s)\right\Vert _{Y}^{4}ds=\frac{M^{4}}{4m^{2}}e^{8m^{2}T}~\mathbb{E}%
\int_{0}^{\theta }\left\Vert Z(s)\right\Vert _{Y}^{4}ds<\infty
\end{multline*}%
because of condition (\ref{chi}), with $M$ and $m$ from (\ref{U-esttt}). The
proof is complete.
\hfill%
$\square $

\begin{remark}
\label{Ito-rem}It can be shown by similar arguments that, if a process $%
y(t)\in Y$, $\ t\in \left[ 0,\theta \right] $, is a solution of integral
equation (\ref{Ito-y}), then
\begin{equation*}
Z(t):=U(w(t))y(t)\in Y
\end{equation*}%
satisfies (\ref{Ito-z}).
\end{remark}

Now we can proceed with the proof of the main result of this section.

\noindent \textbf{Proof of Theorem \ref{ode-theor}. }Let $\theta $ be a
finite stopping time and $z(t),\ t\in \lbrack 0,\theta ]$, a $Y$-valued
process solving the integral equation (\ref{ODE-int}). It is clear that $%
y(t)=U(-w(t))z(t)\in Y$ is a solution of the equation%
\begin{equation}
y(t)=U(-w(t))\left( y_{0}-\int_{0}^{t}U(w(s))F(y(s))ds\right) ,\ t\in
\lbrack 0,\theta ].  \label{int-eq2}
\end{equation}%
We can now apply Lemma \ref{ito-lemma} with $Z(t)=z(t)$, $Z_{0}=y_{0}$ and%
\begin{equation*}
\chi (t)=U(w(t))F(y(t)),\ t\in \lbrack 0,\theta ],
\end{equation*}%
and obtain
\begin{multline*}
y(t\wedge \theta )=\int_{0}^{t\wedge \theta }U(-w(s))\chi
(s)ds-\int_{0}^{t\wedge \theta }DU(-w(s))Z(s)\circ dw(s) \\
=\int_{0}^{t\wedge \theta }F(y(s))ds-\int_{0}^{t\wedge \theta
}DU(-w(s))Z(s)\circ dw(s) \\
=\int_{0}^{t\wedge \theta }F(y(s))ds-\int_{0}^{t\wedge \theta }Dy(s)\circ
dw(s).
\end{multline*}%
The converse implication can be shown by similar arguments, cf. Remark \ref%
{Ito-rem}.
\hfill%
$\square $

\end{document}